\documentclass[12pt]{article}
\usepackage{graphicx}
\usepackage{amsmath}
\usepackage{amsfonts}
\usepackage{amssymb}

\newtheorem{theorem}{Theorem}

\newtheorem{corollary}[theorem]{Corollary}

\newtheorem{lemma}[theorem]{Lemma}

\newtheorem{proposition}[theorem]{Proposition}

\newenvironment{proof}[1][Proof]{\textbf{#1.} }{\ \rule{0.5em}{0.5em}}

\def\a{\alpha}
\def\b{\beta}
\def\g{\gamma}

\def\ps{\varphi}

\def\o{\omega}

\def\S{{\bf S}}

\def\A{{\cal A}}
\def\B{{\cal B}}
\def\A'{{\cal A}'}
\def\B'{{\cal B}'}

\def\Z{{\mathbf Z}}

\date{}

\begin{document}

\title{Seifert manifolds and $(1,1)$-knots}

\author{Luigi Grasselli \and Michele Mulazzani}

\maketitle
\begin{abstract}
{The aim of this paper is to investigate the relations between
Seifert manifolds and $(1,1)$-knots. In particular, we prove that
every orientable Seifert manifold with invariants
$$
\begin{array}{ccc}
\{Oo,0\mid -1;\!\!\!\!\!& \underbrace{(p,q), \ldots, (p,q)},& \!\!\!\!\!(l, l-1)\}\\
&n \mbox{ times}&
\end{array}
$$
has the fundamental group cyclically presented by
$G_n((x_1^q\cdots x_n^q)^lx_n^{-p})$ and, moreover, it is the
$n$-fold strongly-cyclic covering of the lens space $L(\vert
nlq-p\vert,q)$, branched over the $(1,1)$-knot
$K(q,q(nl-2),p-2q,p-q)$, if $p\ge 2q$, and over the $(1,1)$-knot
$K(p-q,2q-p,q(nl-2),p-q)$, if $p<2q$.
}

\end{abstract}

\bigskip

\section{Introduction}

Cyclic branched coverings of knots in $\S^3$ with cyclically
presented fundamental group have been deeply investigated in the
recent years by many authors (see
\cite{BKM,CHK,CHK2,CHR,HKM2,Ki,KKV1,KKV2,MR,VK}). Their results
have been included in an organic and more general context in
\cite{Mu}, where it is proved that the fundamental group of every
$n$-fold strongly-cyclic branched covering of a $(1,1)$-knot
admits a cyclic presentation encoded by a genus $n$ Heegaard
diagram. In \cite{CM1} this result has been improved, obtaining a
constructive algorithm which explicitly gives the cyclic
presentation, starting from a representation of the $(1,1)$-knot
via the mapping class group of the twice punctured torus (see
\cite{CM2} for further details on this representation).

In \cite{Du}, M. J. Dunwoody introduces a class of manifolds
depending on six integers, the so-called Dunwoody manifolds, with
cyclically presented fundamental groups. It is proved in \cite{GM}
and \cite{CM3} that the family of Dunwoody manifolds coincides
with the family of strongly-cyclic coverings of lens spaces
(possibly $\S^3$), branched over $(1,1)$-knots. As a consequence,
any $(1,1)$-knot can be represented by four integers $a,b,c,r$,
and it will be denoted by $K(a,b,c,r)$.

In this paper we show that the orientable Seifert manifold with
invariants
$$
\begin{array}{ccc}
\{Oo,0\mid -1;\!\!\!\!\!& \underbrace{(p,q), \ldots, (p,q)},& \!\!\!\!\!(l, l-1)\}\\
&n \mbox{ times}&
\end{array}
$$
has the fundamental group isomorphic to the cyclically presented
group $G_n((x_1^q\cdots x_n^q)^lx_n^{-p})$, and it is the $n$-fold
strongly-cyclic covering of the lens space $L(\vert
nlq-p\vert,q)$, branched over the $(1,1)$-knot
$K(q,q(nl-2),p-2q,p-q)$, if $p\ge 2q$, and over the $(1,1)$-knot
$K(p-q,2q-p,q(nl-2),p-q)$, if $p< 2q$.

\section{Basic notions}

A finite balanced presentation of a group $<x_1,\ldots,x_n\mid
r_1,\ldots,r_n>$ is said to be a {\it cyclic presentation\/} if
there exists a word $w$ in the free group $F_n$ generated by
$x_1,\ldots,x_n$ such that $r_k=\theta^{k-1}(w)$, $k=1,\ldots,n$,
where \hbox{$\theta :F_n\to F_n$} denotes the automorphism defined
by $\theta (x_i)=x_{i+1}$ (subscripts mod $n$), $i=1,\ldots,n$.
This presentation (and the related group) will be denoted by
$G_n(w)$. For further details see \cite{Jo}.

\begin{figure}[ht]
\begin{center}
\includegraphics*[totalheight=3cm]{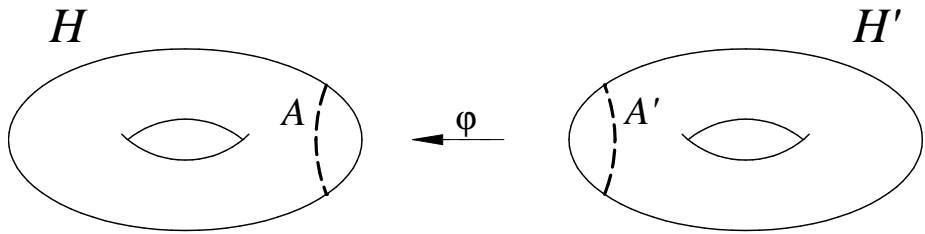}
\end{center}
\caption{A $(1,1)$-knot decomposition.} \label{Fig. 1}
\end{figure}

A knot $K$ in a closed, connected, orientable 3-manifold $N^3$ is
called a $(1,1)$-{\it knot\/} if there exists a Heegaard splitting
of genus one \hbox{$(N^3,K)=(H,A)\cup_{\ps}(H',A'),$} where $H$
and $H'$ are solid tori, $A\subset H$ and $A'\subset H'$ are
properly embedded trivial arcs\footnote{This means that there
exists a disk $D\subset H$ (resp. $D'\subset H'$) with $A\cap
D=A\cap\partial D=A$ and $\partial D-A\subset\partial H$ (resp.
$A'\cap D'=A'\cap\partial D'=A'$ and $\partial
D'-A'\subset\partial H'$).}, and $\ps:(\partial H',\partial
A')\to(\partial H,\partial A)$ is an attaching homeomorphism (see
Figure \ref{Fig. 1}). Obviously, $N^3$ turns out to be a lens
space $L(p,q)$, including $\S^3=L(1,0)$ and
$\S^2\times\S^1=L(0,1)$.

It is well known that the family of $(1,1)$-knots contains all
torus knots and all two-bridge knots in $\S^3$. Several
topological properties of $(1,1)$-knots have recently been
investigated in many papers (see references in \cite{CM2}).

An integer 4-parametric representation of $(1,1)$-knots have been
developed in \cite{CM3} (see also \cite{CM4}). Every $(1,1)$-knot,
with the only exception of the ``core'' knot
$\{P\}\times\S^1\subset\S^2\times\S^1$, can be represented by four
non-negative integers $a,b,c,r$, and the represented knot will be
referred as $K(a,b,c,r)$.

\begin{figure}
 \begin{center}
 \includegraphics*[totalheight=9cm]{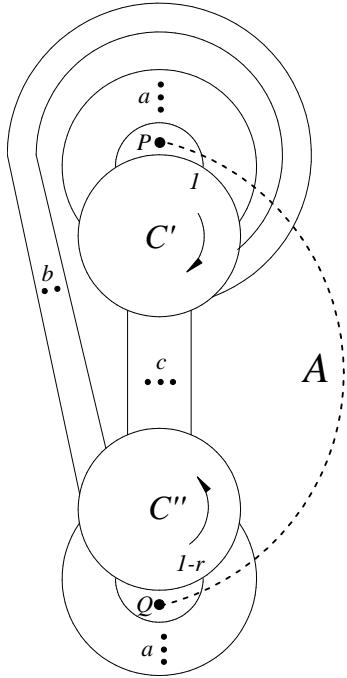}
 \end{center}
 \caption{Heegaard diagram for $K(a,b,c,r)$}
 \label{parametrization}
\end{figure}

A $(1,1)$-knot $K(a,b,c,r)$, with $a+b+c>0$, admits a natural
$(1,1)$-decomposition $(H,A)\cup_{\ps}(H',A')$ described by the
genus one Heegaard diagram of Figure \ref{parametrization}, where
the labels $a,b$ and $c$ denote the correspondent number of
parallel arcs, and the gluing between the circles $C'$ and $C''$
depends on the twist parameter $r$ in such a way that equally
labelled vertices are identified together (observe that $r$ can be
taken mod $2a+b+c$).

\begin{figure}
 \begin{center}
 \includegraphics*[totalheight=3cm]{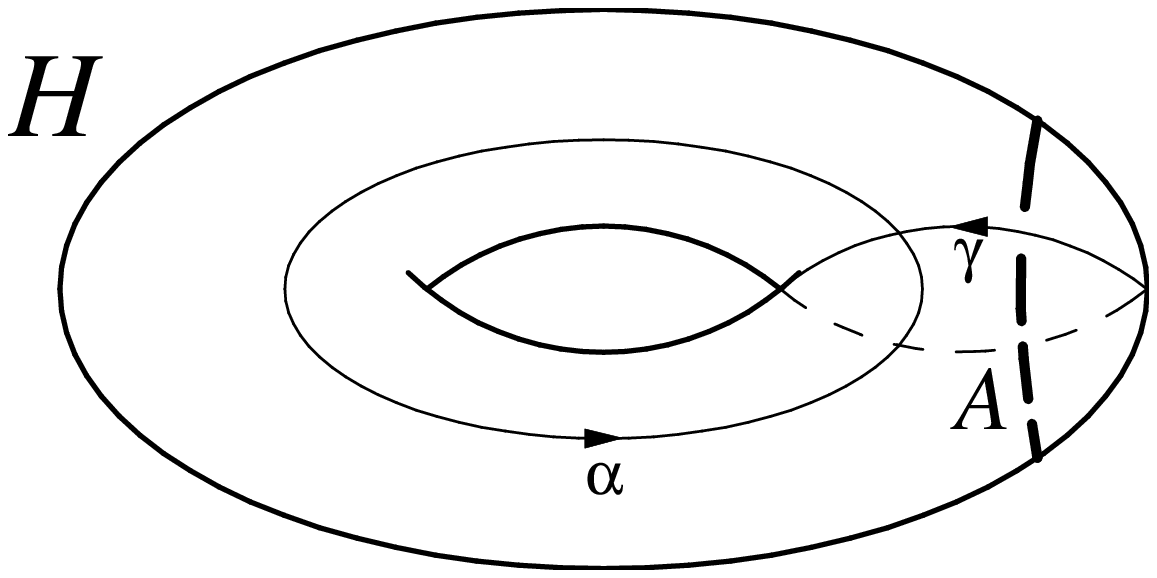}
 \end{center}
 \caption{}
 \label{generators}
\end{figure}

As a simple consequence of Seifert-Van Kampen Theorem, the
fundamental group of the exterior of a $(1,1)$-knot (as well as
its first homology group) is generated by the two loops
$\a,\g\subset\partial H$ depicted in Figure \ref{generators}.

In the next section we need the following result.

\begin{lemma} \label{lens}
\begin{enumerate}

\item[(i)] If $a$ and $c$ are non-negative integers such that
\hbox{$\gcd(a,c)=1$,} then $K(a,0,c,a)$ is a $(1,1)$-knot in the
lens space $L(c,a)$.

\item[(ii)] If $a,b,c$ are non-negative integers such that
\hbox{$\gcd(a-c,b+c)=1$,} then $K(a,b,c,a)$ is a $(1,1)$-knot in
the lens space $L(b+c,a+b)$.

\item[(iii)] If $a,b,c$ are non-negative integers such that $a>0$
and \hbox{$\gcd(a,b-c)=1$}, then $K(a,b,c,a+c)$ is a $(1,1)$-knot
in the lens space $L(\vert b-c \vert,a)$.

\end{enumerate}
\end{lemma}

\begin{figure}
 \begin{center}
 \includegraphics*[totalheight=8cm]{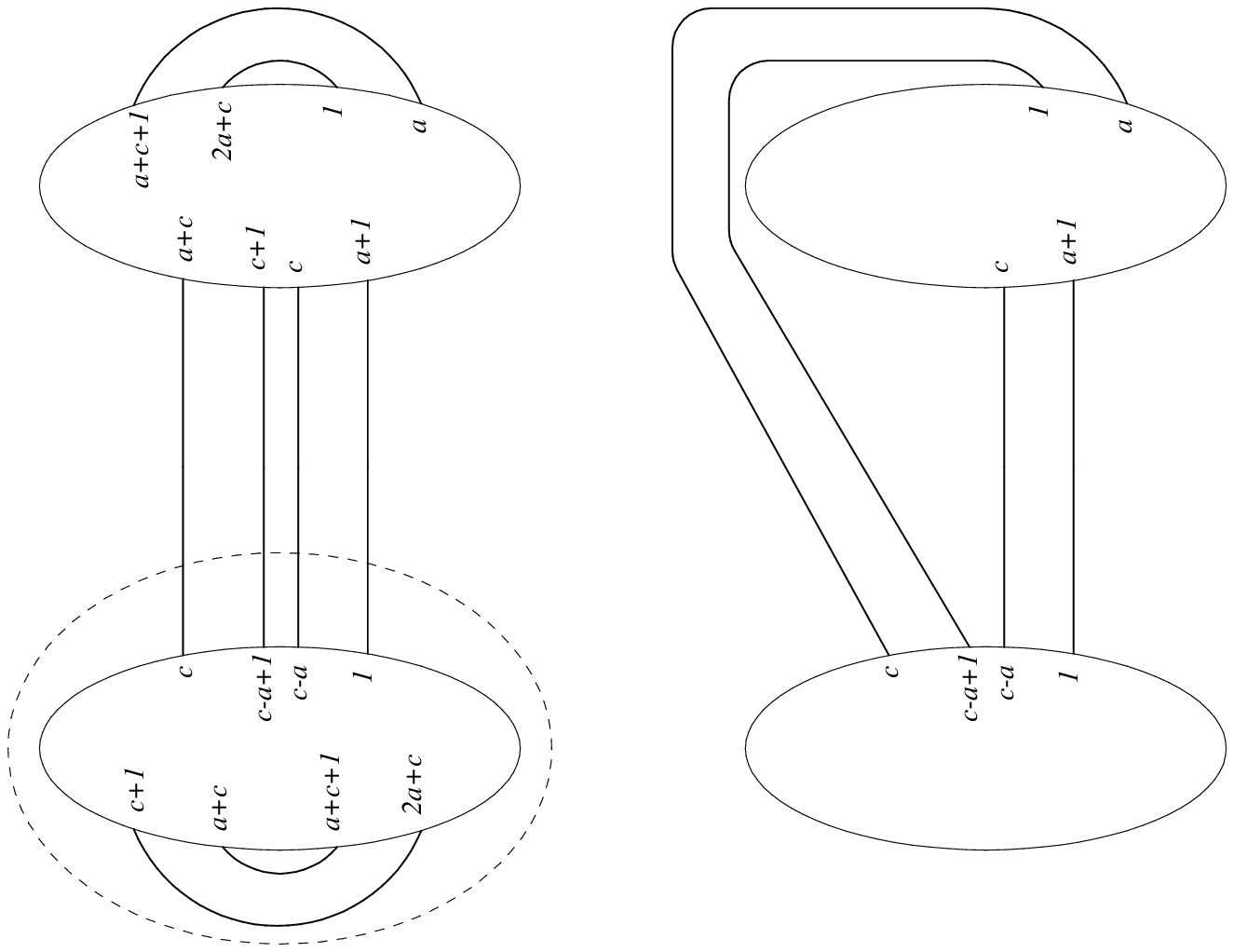}
 \end{center}
 \caption{}
 \label{Singer2}
\end{figure}

\begin{figure}
 \begin{center}
 \includegraphics*[totalheight=8cm]{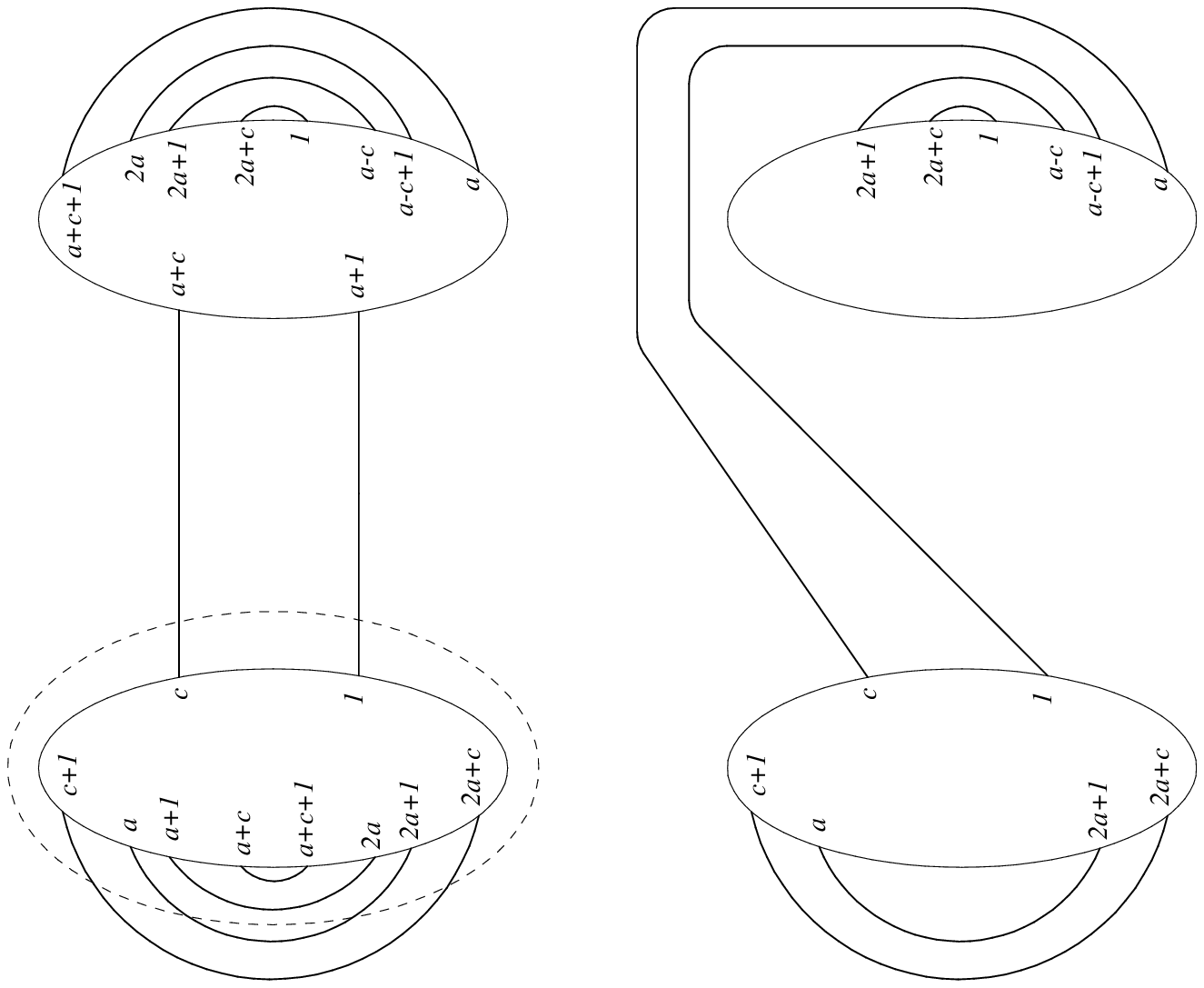}
 \end{center}
 \caption{}
 \label{Singer1}
\end{figure}

\begin{figure}
 \begin{center}
 \includegraphics*[totalheight=8cm]{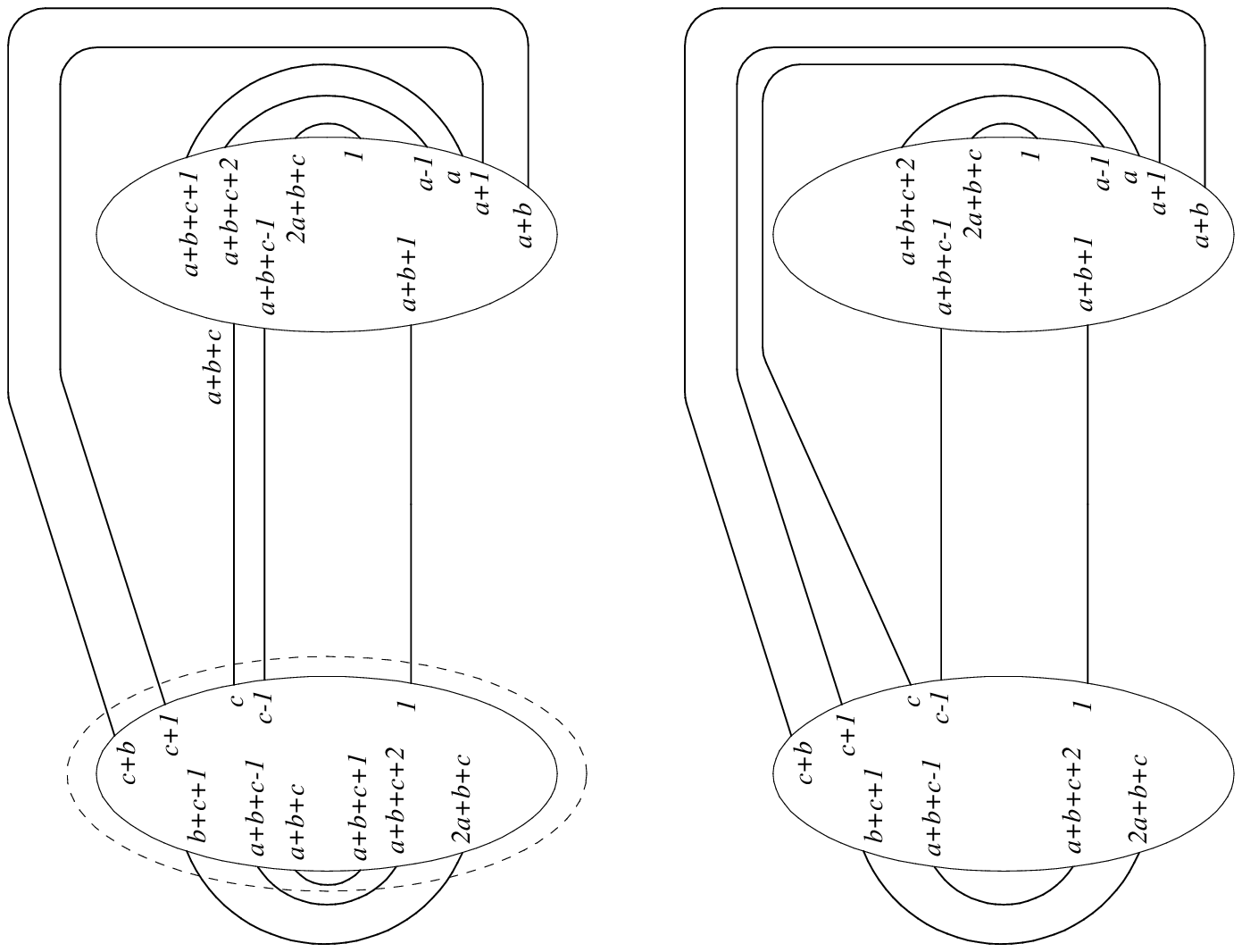}
 \end{center}
 \caption{}
 \label{Singer3}
\end{figure}

\begin{figure}
 \begin{center}
 \includegraphics*[totalheight=8cm]{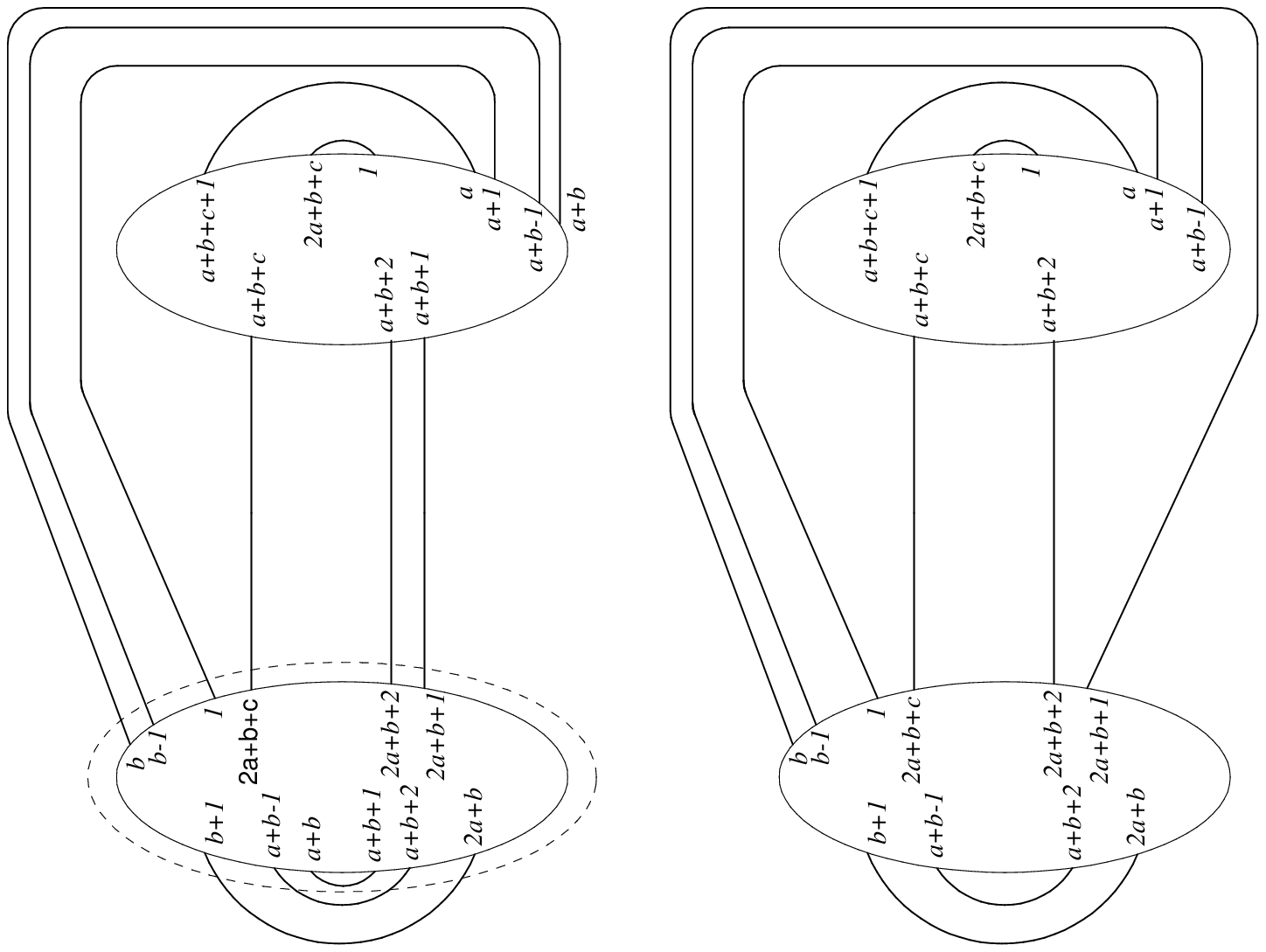}
 \end{center}
 \caption{}
 \label{Singer4}
\end{figure}

\begin{proof} As proved in \cite{CM4}, $K(a,b,c,r)$ is
equivalent to $K(a,c,b,2a+b+c-r)$ and $K(a,0,c,r)$ is equivalent
to $K(a,c,0,r)$. As a consequence, $K(a,0,c,r)$ is equivalent to
$K(a,0,c,2a+c-r)$.

(i) If $a=0$, then $c=1$ and the result is straightforward. If
$0<a<c$, then applying the Singer move \cite{Si} of type IB
depicted in Figure \ref{Singer2}, we obtain the canonical Heegaard
diagram of $L(c,a)$. If $a\ge c$, the Singer move depicted in
Figure \ref{Singer1} transforms the diagram of $K(a,0,c,a)$ in the
diagram of $K(a-c,c,0,a-c)$, which is equivalent to
$K(a-c,0,c,a-c)$. Since $L(c,a-c)$ is homeomorphic to $L(c,a)$,
the result follows from the previous case $a<c$.

(ii) If $b=0$ the result follows from (i). If $b>0$, by performing
the Singer move of Figure \ref{Singer3}, the diagram of
$K(a,b,c,a)$ becomes the diagram of $K(a-1,b+1,c-1,a-1)$. If $a\le
c$, after performing the move $a$ times, we obtain the diagram of
$K(0,b+a,c-a,0)$, which is the canonical Heegaard diagram of
$L(b+c,a+b)$, since $\gcd(a+b,b+c)=\gcd(a-c,b+c)=1$. If $a>c$,
after performing the move $c$ times, we obtain the diagram of
$K(a-c,b+c,0,a-c)$, which is equivalent to $K(a-c,0,b+c,a-c)$. Now
the result follows from (i), since $L(b+c,a-c)$ is homeomorphic to
$L(b+c,a+b)$.

(iii) Since $K(a,b,c,a+c)$ is equivalent to $K(a,c,b,a+b)$, we can
always suppose $c\le b$. If $c>0$, then, by performing the Singer
move of Figure \ref{Singer4}, the diagram of $K(a,b,c,a+c)$
becomes the diagram of $K(a,b-1,c-1,a+c-1)$. After performing the
move $c$ times, we obtain the diagram of $K(a,b-c,0,a)$, which is
equivalent to $K(a,0,b-c,a)$. The result now follows from (i).
\end{proof}

\medskip

Observe that, since $K(a,b,c,a+b+c)$ is equivalent to
$K(a,c,b,a)$, $K(a,b,c,a+b+c)$ is a $(1,1)$-knot in the lens space
$L(b+c,a+c)$, if \hbox{$\gcd(a-b,b+c)=1$} .

An $n$-fold cyclic covering $M^3$ of a 3-manifold $N^3$, branched
over a knot $K\subset N^3$, is called {\it strongly-cyclic\/} if
the branching index of $K$ is $n$. This means that the fiber in
$M^3$ of each point of $K$ consists of a single point. In this
case, the homology class $m$ of a meridian loop around $K$ is
mapped by the monodromy $\o:H_1(N^3-K)\to \Z_n$ of the covering in
a generator of $\Z_n$ (up to equivalence, we can always suppose
$\o(m)=1$). Observe that a cyclic branched covering of a knot $K$
in $\S^3$ is always strongly-cyclic and is uniquely determined, up
to equivalence, since $H_1(\S^3-K)\cong\Z$. Obviously, this
property is no longer true for a knot in a more general
3-manifold. Necessary and sufficient conditions for the existence
and uniqueness of strongly-cyclic branched coverings of
$(1,1)$-knots have been obtained in \cite{CM1}.

\section{Main results}

Let $n,p,q,l$ be positive integers such that $q<p$ and
$\gcd(p,q)=1$. In the following we denote by $\Sigma(n,p,q,l)$ the
orientable Seifert manifold \cite{ST} with invariants
$$
\begin{array}{ccc}
\{Oo,0\mid -1;\!\!\!\!\!& \underbrace{(p,q), \ldots, (p,q)},& \!\!\!\!\!(l, l-1)\}\\
&n \mbox{ times}&
\end{array}
$$
having $\S^2$ as orbit space, $n$ exceptional fibres of type
$(p,q)$ and, for $l>1$, an exceptional fibre of type $(l,l-1)$.

Observe that, in particular, the manifolds $\Sigma(n,2,1,1)$ are
precisely the Neuwirth manifolds $M_n$ introduced in \cite{Ne},
studied in \cite{Ca}, and generalized in \cite{GP}, \cite{SV} and
\cite{RSV}.

\begin{proposition} \label{fundamental}
The fundamental group of $\Sigma(n,p,q,l)$ is cyclically presented
by $G_n(w)$, where $w=(x_1^q\cdots x_n^q)^lx_n^{-p}$.
\end{proposition}

\begin{proof}
Following \cite{Or}, a standard presentation of
$G=\pi_1(\Sigma(n,p,q,l))$ is
$$\langle y_1,\ldots,y_n,y,h\mid [y_i,h],[y,h],y_i^ph^q,y^lh^{l-1},y_1\cdots y_nyh; i=1,\ldots,n\rangle.$$
Since $\gcd(p,q)=1$, there exist $\a,\b\in\Z$ such that
$q\b-p\a=1$.

From the relations we have $y^lh^{l-1}=(yh)^lh^{-1}$. By
introducing the new generator $x=yh$, the presentation becomes
$$\langle y_1,\ldots,y_n,x,h\mid
[y_i,h],[x,h],y_i^ph^q,x^lh^{-1},y_1\cdots y_nx;
i=1,\ldots,n\rangle.$$

Now we define $x_i=y_i^{\b}h^{\a}$, for $i=1,\ldots,n$. Then we
have the following new presentation for $G$:
$$\langle y_1,\ldots,y_n,x,h,x_1,\ldots,x_n\mid
[y_i,h],[x,h],y_i^ph^q,x^lh^{-1},y_1\cdots
y_nx,x_i^{-1}y_i^{\b}h^{\a}; i=1,\ldots,n\rangle.$$

From the relations we obtain
$x_i^q=y_i^{q\b}h^{q\a}=y_i^{1+p\a}h^{q\a}=y_i(y_i^ph^q)^{\a}=y_i$
and
$x_i^p=y_i^{p\b}h^{p\a}=y_i^{p\b}h^{q\b-1}=(y_i^ph^q)^{\b}h^{-1}=h^{-1}.$

Therefore $G$ admits the presentation
$$\langle x,h,x_1,\ldots,x_n\mid
[x_i^q,x_i^{-p}],[x,x^l],x_i^{qp}x_i^{-qp},x^lh^{-1},x_1^q\cdots
x_n^qx,x_i^{-1}x_i^{q\b}x_i^{-p\a},x_i^ph; i=1,\ldots,n\rangle,$$
which obviously becomes
$$\langle x,h,x_1,\ldots,x_n\mid
x^lh^{-1},x_1^q\cdots x_n^qx,x_i^ph; i=1,\ldots,n\rangle.$$

From the relations we obtain $x=(x_1^q\cdots x_n^q)^{-1}$ and
$h=(x_1^q\cdots x_n^q)^{-l}$.

Therefore $G$ admits the presentation:
$$\langle x_1,\ldots,x_n\mid
(x_1^q\cdots x_n^q)^lx_i^{-p}; i=1,\ldots,n\rangle,$$ which is
equivalent to
$$\langle x_1,\ldots,x_n\mid
(x_1^q\cdots x_n^q)^lx_n^{-p},x_i^px_{i+1}^{-p};
i=1,\ldots,n-1\rangle.$$

In order to prove that $G$ is isomorphic to $G_n(w)$, it suffices
to show that the normal closures in the free group $F_n$ of the
sets $\{w,x_i^px_{i+1}^{-p}; i=1,\ldots,n-1\}$ and
$\{w,\theta^i(w);i=1,\ldots,n-1\}$ coincide.

This is achieved by observing that, setting $w_i=\theta^i(w)$, we
have:

$w_i^{-1}x_{i+1}^qw_{i+1}x_{i+1}^{-q}=x_i^px_{i+1}^{-p}$, for
$i=1,\ldots,n-2$,

$w_{n-1}^{-1}x_n^qwx_n^{-q}=x_{n-1}^px_n^{-p}$,

$x_1^{-q}w((x_1^px_2^{-p})\cdots(x_{n-1}^px_n^{-p}))^{-1}x_1^q=w_1$,

$x_i^{-q}w_{i-1}(x_{i-1}^px_i^{-p})x_i^q=w_i$, for
$i=2,\ldots,n-1$.

\end{proof}

Since a Seifert manifold with one or two singular fibers is
trivially a lens space (see \cite{Or}), we will suppose $n>1$ and
$l>1$ when $n=2$.

\begin{theorem} \label{general}
The Seifert manifold $\Sigma(n,p,q,l)$ is the $n$-fold
strongly-cyclic covering of the lens space $L(\vert
nlq-p\vert,q)$, branched over the $(1,1)$-knot
$K(q,q(nl-2),p-2q,p-q)$, if $p\ge 2q$, and over the $(1,1)$-knot
\hbox{$K(p-q,2q-p,q(nl-2),p-q)$,} if $p< 2q$.
\end{theorem}

\begin{figure}
 \begin{center}
 \includegraphics*[totalheight=6.5cm]{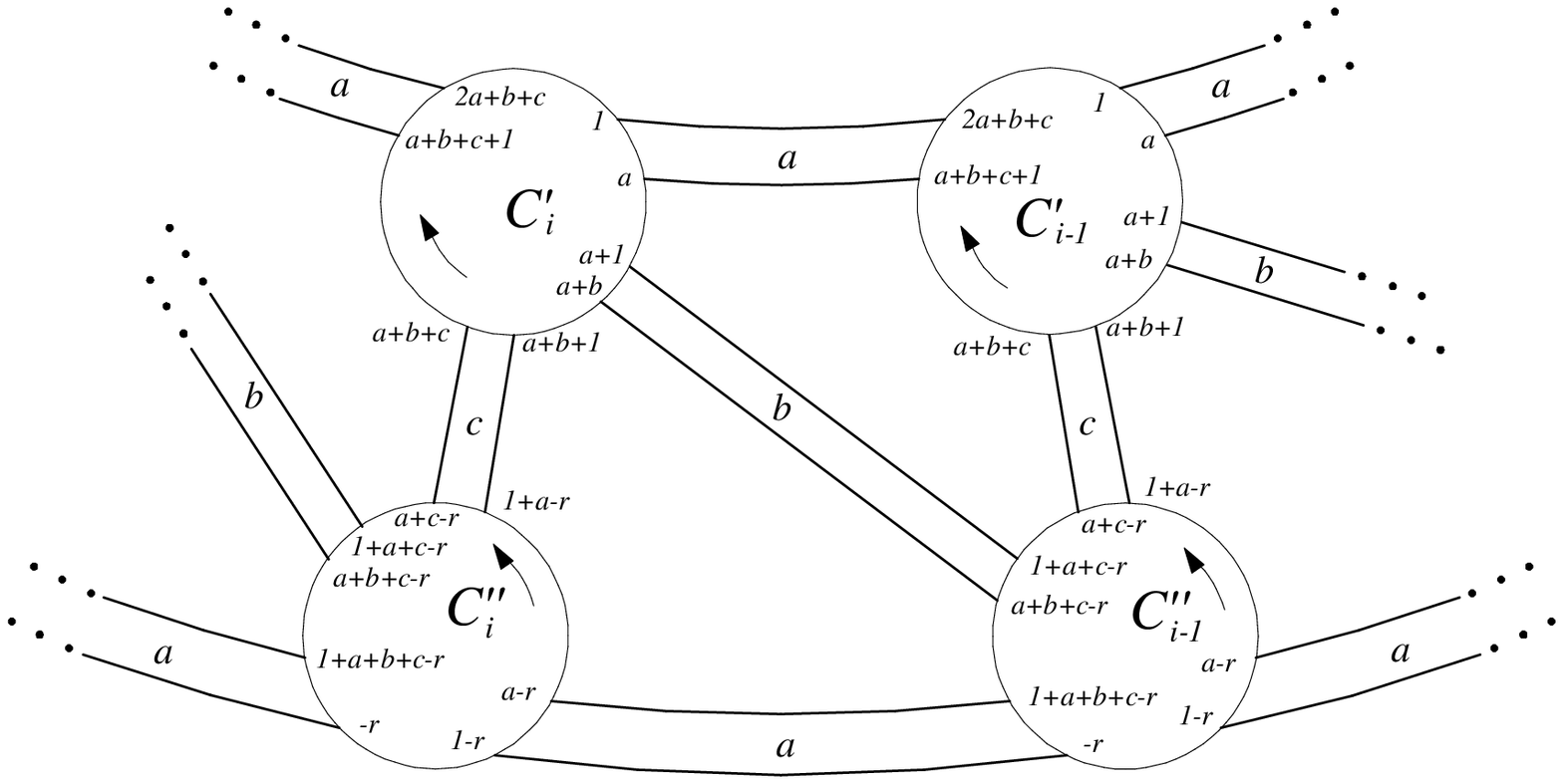}
 \end{center}
 \caption{}
 \label{Dunwoody}
\end{figure}

\begin{figure}
 \begin{center}
 \includegraphics*[totalheight=8cm]{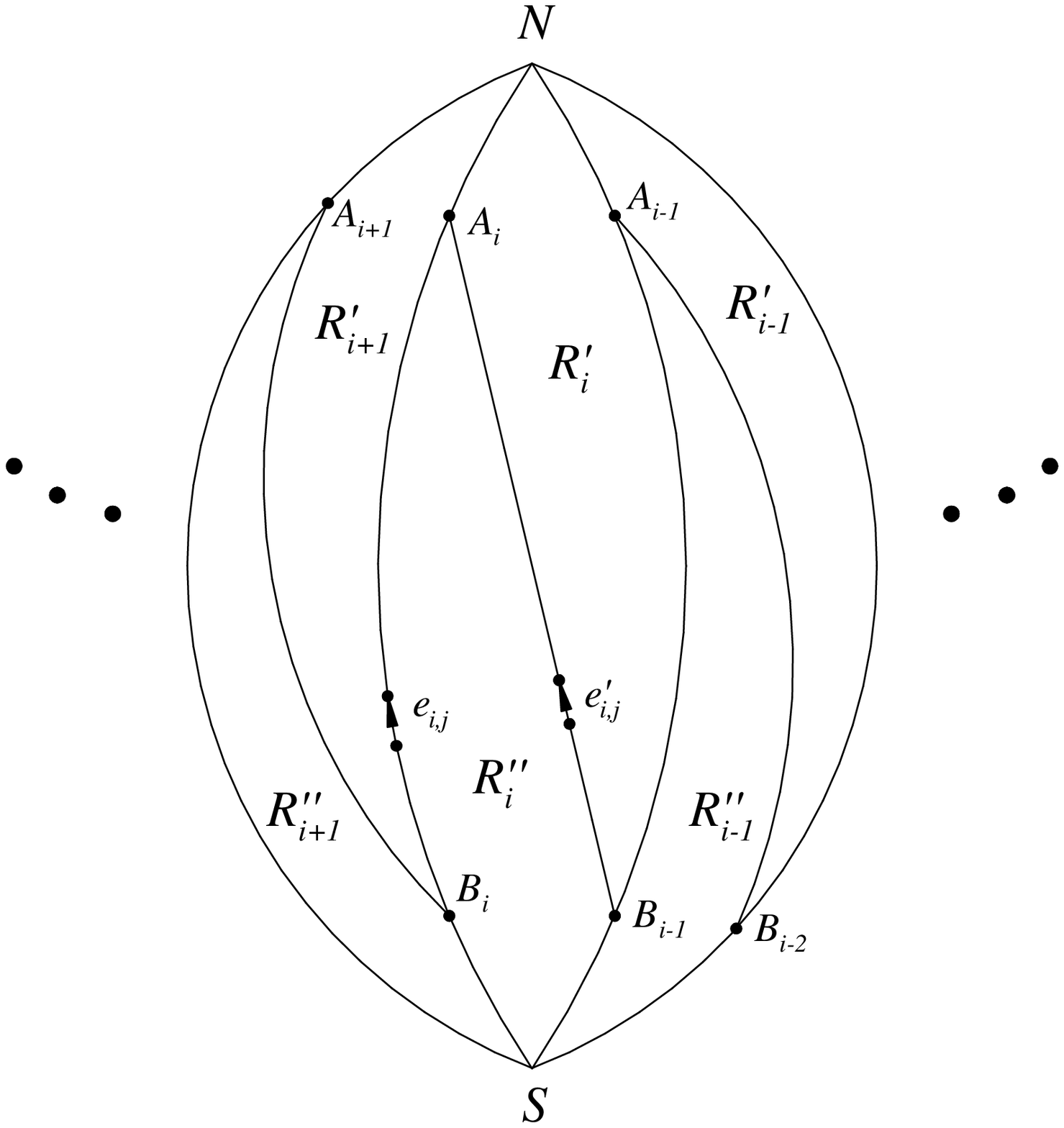}
 \end{center}
 \caption{}
 \label{tess}
\end{figure}

\begin{figure}
 \begin{center}
 \includegraphics*[totalheight=8cm]{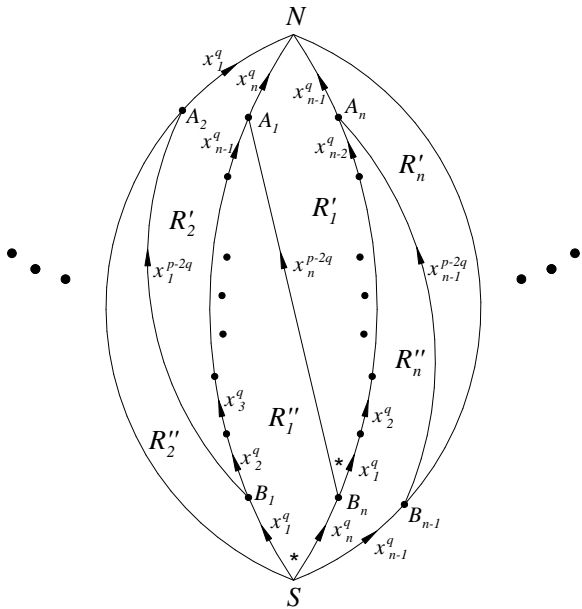}
 \end{center}
 \caption{Case $p\ge 2q$}
 \label{tess1}
\end{figure}

\begin{figure}
 \begin{center}
 \includegraphics*[totalheight=8cm]{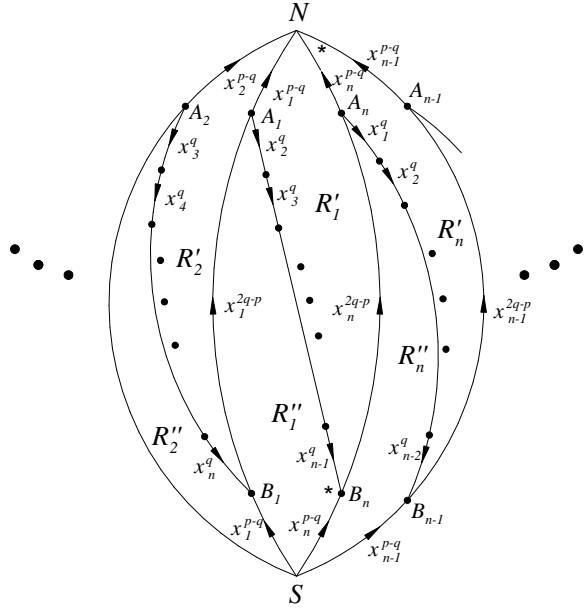}
 \end{center}
 \caption{Case $p < 2q$}
 \label{tess2}
\end{figure}

\begin{proof}
(i) Suppose $p\ge 2q$. By (iii) of Lemma \ref{lens},
$K=K(q,q(nl-2),p-2q,p-q)$ is a $(1,1)$-knot in $L(\vert
nlq-p\vert,q)$. Now suppose that $K$ admits the $n$-fold
strongly-cyclic branched covering $C_n(K)$ defined by $\o(\a)=0$
and $\o(\g)=1$. In this case $C_n(K)$ is the Dunwoody manifold
$D(q,q(nl-2),p-2q,n,p-q,0)$ (see \cite{GM} and \cite{CM3}). Its
defining genus $n$ Heegaard diagram, having cyclic symmetry of
order $n$, is depicted in Figure \ref{Dunwoody}, where
$a=q,b=q(nl-2),c=p-2q$, $r=p-q$, and the circle $C_i'$ must be
glued to the circle $C_i''$, for $i=1,\ldots,n$, according to the
twist parameter $r$. In order to check if the above diagram really
represents a manifold, it is convenient to consider the cellular
decomposition dual to the one associated to the diagram. In this
way, $C_n(K)$ is obtained by pairwise identification of the
regions of the 2-cell tessellation of the boundary of a 3-ball, as
depicted in Figure \ref{tess}. The tessellation consists of $2n$
regions, $R'_1,\ldots,R'_n$ around the North pole $N$ and
$R''_1,\ldots,R''_n$ around the South pole $S$, corresponding to
the cycles $C_i'$ and $C_i''$ of the diagram. The edges of the
tessellation correspond to the arcs of the diagram in the
following way. The poles $S$ and $N$ are connected by $n$
meridians and each meridian $m_i$ is composed by $nlq$ edges
$e_{i,j}$ ($j=1,\ldots, nlq$) from $S$ to $N$, and precisely $q$
from $S$ to $B_i$, $q(nl-2)$ from $B_i$ to $A_i$ and $q$ from
$A_i$ to $N$. Moreover, for $i=1,\ldots,n$, the arc $a_i$
connecting $B_{i-1}$ with $A_i$ is composed by $p-2q$ edges
$e'_{i,j}$ ($j=1,\ldots, p-2q$). Observe that, when $p=2q$ (i.e.,
$p=2$ and $q=1$), the points $B_{i-1}$ and $A_i$ coincide. In
order to obtain $C_n(K)$, the region $R'_i$ is glued to $R''_i$,
for $i=1,\ldots,n$, by an orientation reversing homeomorphism in
such a way that point $S$ of $R''_i$ matches point\footnote{The
subscript $i$ is considered mod $n$.} $B_{i-1}$ of $R'_i$. By this
gluing we have: $e_{i,j}\equiv e_{i-1,j+q}$ for
$j=1,\ldots,q(nl-1)$, $e_{i-1,j}\equiv e'_{i,j}$ for
$j=1,\ldots,q$, $e'_{i,j}\equiv e'_{i,j+q}$ for $j=1,\ldots,p-3q$,
and $e'_{i,p-3q+k}\equiv e_{i,q(nl-1)+k}$ for $k=1,\ldots,q$. As a
consequence of all identifications, we have $e_{i,j}\equiv
e_{i,j+nq}\equiv \cdots \equiv e_{i,j+n(l-1)q}$ for
$j=1,\ldots,nq$, and moreover $e_{i-1,j}\equiv e_{i,q(nl-1)+j}$
for $j=1,\ldots,q$. Since $\gcd(p,q)=1$, all edges of the arcs
$SB_{i-1}$, $B_{i-1}A_i$ and $A_iN$ match, and we will call them
$x_i$. Therefore, the $q(nl-2)$ edges of the arc $B_iA_i$ become:
$q$ times $x_{i+1}$, $q$ times $x_{i+2}$, \ldots, $q$ times
$x_{i}$ ($l-1$ times) and then $q$ times $x_{i+1}$, $q$ times
$x_{i+2}$, \ldots, $q$ times $x_{i-2}$ (see Figure \ref{tess1}).
Since each edge $x_i$ appears $p$ times consecutively, its
endpoints coincide. Moreover, since every edge $x_i$ has $S$ as
endpoint, the cellular decomposition of $C_n(K)$ is composed by
one vertex, $n$ edges, $n$ regions and one 3-ball. By Seifert
criterion, $C_n(K)$ is actually a closed, orientable 3-manifold.

Moreover, $\pi_1(C_n(K))$ admits a balanced presentation with
$x_1,\ldots,x_n$ as generators and with relators obtained by
walking around the boundaries of the regions $R'_i$. As a
consequence, $\pi_1(C_n(K))=G_n((x_1^q\cdots x_n^q)^lx_n^{-p})$,
and therefore it is isomorphic to $\pi_1(\Sigma(n,p,q,l))$.

In order to prove that $C_n(K)$ and $\Sigma(n,p,q,l)$ are actually
homeomorphic, we note that $\pi_1(C_n(K))$ has non-trivial centre.
So, either $C_n(K)$ is prime or is the connected sum $C_n(K)=M\#
M'$, where $\pi_1(M')$ is trivial and $M$ is prime, with
$\pi_1(M)=\pi_1(C_n(K))$.

If $\Sigma(n,p,q,l)$ is large (in the sense of \cite{Or}), then
also $M$ is a large Seifert manifold (see \cite{CJ}). So $M$ and
$\Sigma(n,p,q,l)$ are homeomorphic, with Heegaard genus at least
$n-2$ (see \cite{BZ}). If $M'$ is not homeomorphic to $\S^3$, then
its Heegaard genus is greater than two, and $C_n(K)$ has Heegaard
genus greater than $n$. But this is impossible since $C_n(K)$
admits a genus $n$ Heegaard splitting by \cite{Mu}. As a
consequence, $M'=\S^3$ and $C_n(K)=\Sigma(n,p,q,l)$.

If $\Sigma(n,p,q,l)$ is not large, then one of the following
possibilities holds: (a) $n=3$, $p=2$, $q=1=l$; (b) $n=2$, $p=2$,
$q=1$, $l>1$; (c) $n=2$, $p=3$, $q=1$, $l=2$. Since for $q=1$ the
given decomposition of $C_n(K)$ coincides with the decomposition
$P(p,\ldots,p;l)$ given in \cite{RSV}, the result follows from
\cite[Proposition 4.1]{RSV}.

(ii) Let $p<2q$. By (ii) of Lemma \ref{lens},
$K(p-q,2q-p,q(nl-2),p-q)$ is a $(1,1)$-knot in $L(\vert
nlq-p\vert,q)$. Now suppose that $K$ admits the $n$-fold
strongly-cyclic branched covering $C_n(K)$ defined by $\o(\a)=1$
and $\o(\g)=1$. In this case $C_n(K)$ is the Dunwoody manifold
$D(p-q,2q-p,q(nl-2),n,p-q,1)$. Its defining genus $n$ Heegaard
diagram, having cyclic symmetry of order $n$, is depicted in
Figure \ref{Dunwoody}, where $a=p-q,b=2q-p,c=q(nl-2)$, $r=p-q$,
and the circle $C_i'$ must be glued to the circle $C_{i+1}''$, for
$i=1,\ldots,n$, according to the twist parameter $r$. Again it is
better to refer to the dual decomposition (see Figure \ref{tess}).
This time each meridian is composed by $p$ edges, and precisely
$p-q$ from $S$ to $B_i$, $2q-p$ from $B_i$ to $A_i$ and $p-q$ from
$A_i$ to $N$. Moreover, the arc connecting $B_{i-1}$ with $A_i$ is
composed by $q(nl-2)$ edges. In order to obtain $C_n(K)$, the
region $R'_{i-1}$ is glued to $R''_i$ in such a way that point $N$
of $R'_{i-1}$ matches point $B_{i-1}$ of $R''_i$.

\begin{figure}
 \begin{center}
 \includegraphics*[totalheight=6cm]{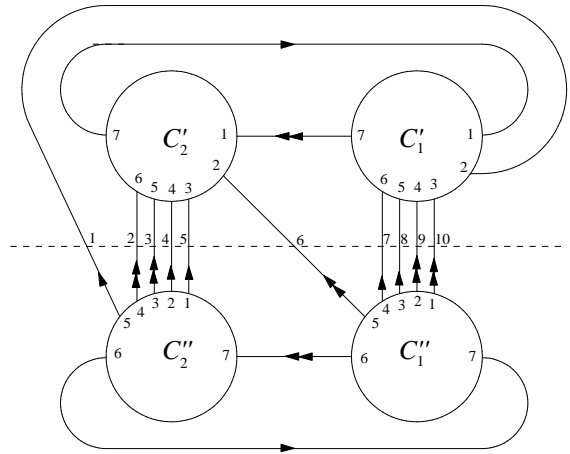}
 \end{center}
 \caption{$D(1,1,4,2,1,1)$}
 \label{knot1}
\end{figure}

\begin{figure}
 \begin{center}
 \includegraphics*[totalheight=6cm]{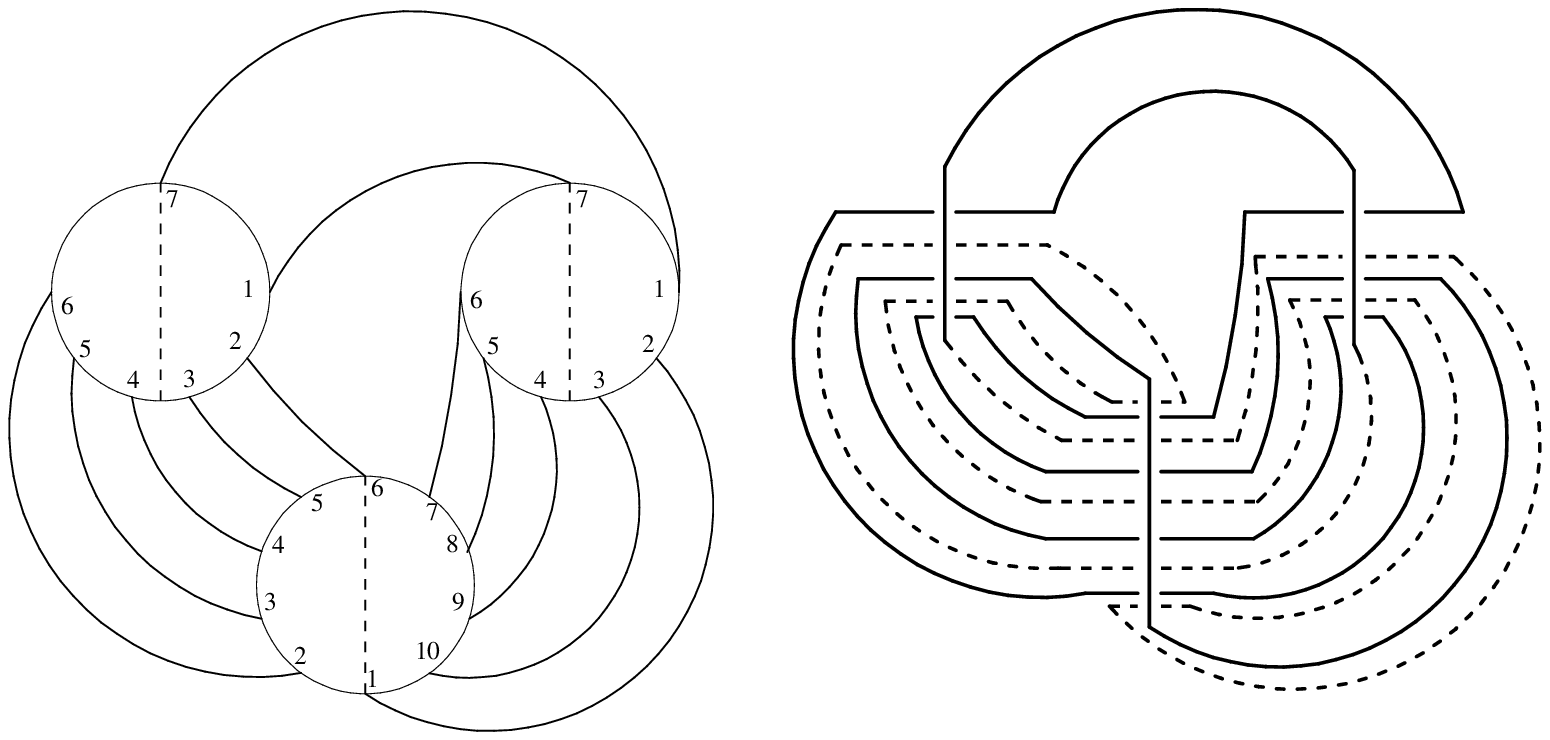}
 \end{center}
 \caption{}
 \label{knot2}
\end{figure}

\begin{figure}
 \begin{center}
 \includegraphics*[totalheight=6cm]{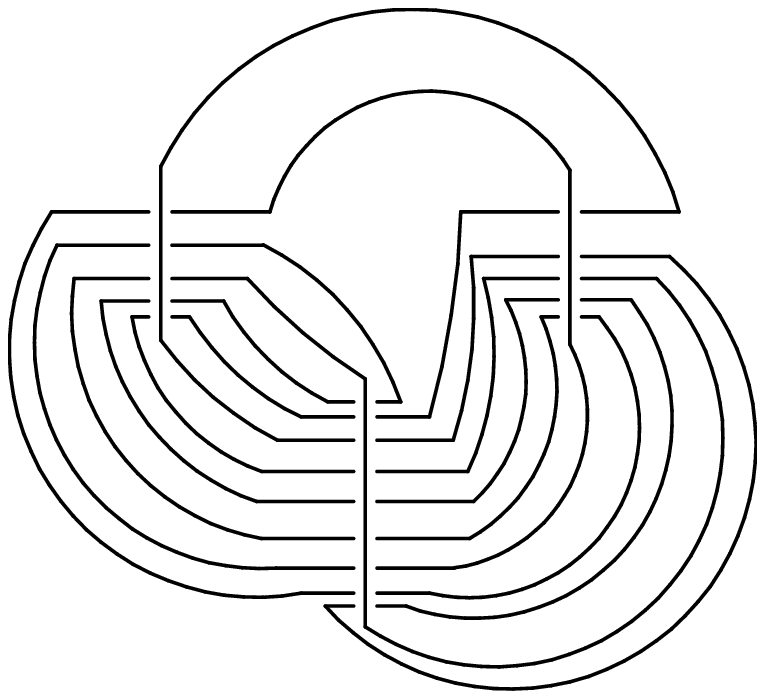}
 \end{center}
 \caption{}
 \label{knot3}
\end{figure}

\begin{figure}
 \begin{center}
 \includegraphics*[totalheight=6cm]{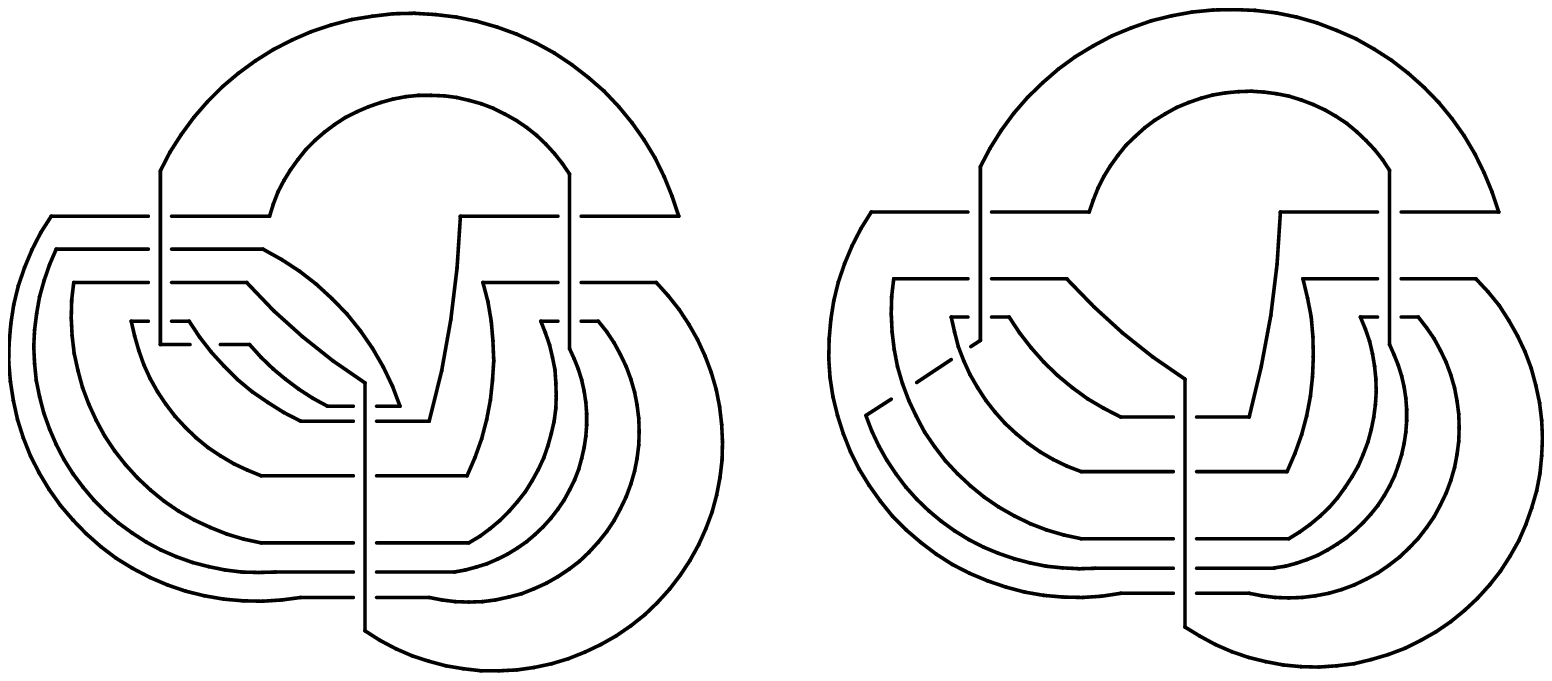}
 \end{center}
 \caption{}
 \label{knot4}
\end{figure}

\begin{figure}
 \begin{center}
 \includegraphics*[totalheight=6cm]{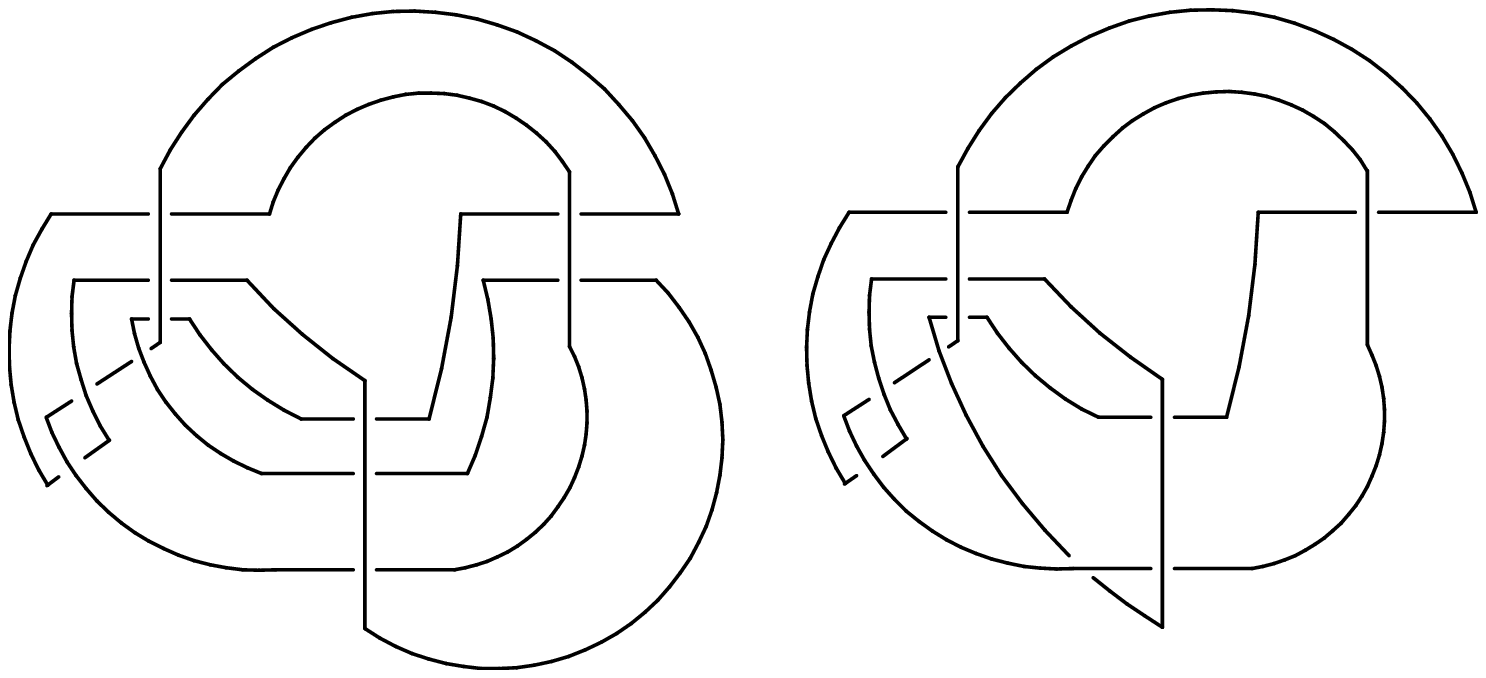}
 \end{center}
 \caption{}
 \label{knot5}
\end{figure}

\begin{figure}
 \begin{center}
 \includegraphics*[totalheight=6cm]{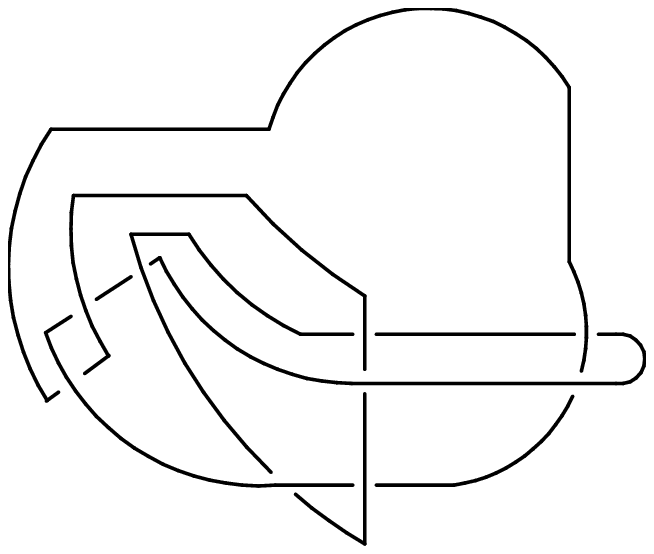}
 \end{center}
 \caption{}
 \label{knot6}
\end{figure}

In complete analogy with point (i), it is easy to see that the
gluing gives rise to a cellular decomposition of $C_n(K)$,
composed by one vertex, $n$ edges, $n$ regions and one 3-ball (see
the 2-skeleton in Figure \ref{tess2}). So $\pi_1(\Sigma(n,p,q,l))$
is isomorphic to $\pi_1(C_n(K))$ and $\Sigma(n,p,q,l)$  is
homeomorphic to $C_n(K)$ when it is large. The only small Seifert
manifold occurs when $n=q=l=2$, $p=3$. To achieve the result, it
suffices to prove that $D=D(1,1,4,2,1,1)$ is the 2-fold covering
of $\S^3$, branched over the Montesinos knot ${\bf
m}(-1;1/2,2/3,2/3)$ (see \cite[Chapter 12]{BuZ}). A genus two
Heegaard diagram of $D$ is depicted in Figure \ref{knot1}, where
the circle $C'_1$ (resp. $C'_2$) must be glued to the circle
$C''_2$ (resp. $C''_1$) in such a way that equally labelled
vertices match. The application of the Takahashi algorithm
\cite{Ta} to this diagram (as depicted in Figure \ref{knot2})
shows that the manifold $D$ is the 2-fold branched covering of the
knot $K$ represented by the 3-bridge diagram of Figure
\ref{knot3}. A Wirtinger presentation of the fundamental group of
the exterior of $K$ can be easily computed from this diagram,
obtaining $\pi_1(\S^3-K)=\langle x,y,z \mid
yz^{-1}xzy^{-1}xyz^{-1}x^{-1}zy^{-1}z^{-1},x^{-1}zyz^{-1}xzx^{-1}zy^{-1}z^{-1}xy^{-1}
\rangle$. The Alexander polynomial of $K$ can be obtained by a
standard application of Fox differential calculus (see
\cite[Chapter 9]{BuZ}), and we have
$\Delta_K(t)=1-4t+5t^2-4t^3+t^4$. By the sequence of Reidemeister
moves depicted in Figures \ref{knot4}-\ref{knot6}, the knot $K$
admits the 9 crossing diagram of Figure \ref{knot6}. Since the
knot $8_{21}$ of Rolfsen's table is the only knot with crossing
number $\le 9$ with Alexander polynomial $1-4t+5t^2-4t^3+t^4$ (see
\cite[Table 3]{Ka}), then $K$ is precisely the knot $8_{21}$.
Since $8_{21}$ is the Montesinos knot ${\bf m}(-1;1/2,2/3,2/3)$
(see \cite[Table 2]{Ka}, which uses a slightly different
notation), the statement is proved.
\end{proof}

Observe that Theorem \ref{general} refines and extends the results
of \cite[Theorem 4.2]{CRS} and \cite[Theorem 5.1]{Sp}. The proof
of Theorem \ref{general} also gives the following:

\begin{corollary}
The Seifert manifold $\Sigma(n,p,q,l)$ is the Dunwoody manifold
$D(q,q(nl-2),p-2q,n,p-q,0)$ if $p\ge 2q$ and the Dunwoody manifold
$D(p-q,2q-p,q(nl-2),n,p-q,1)$ if $p < 2q$.
\end{corollary}


\noindent{\bf Remark} Since
$\Sigma(n,p,p-1,1)=\Sigma(n-1,p,p-1,p)$, this manifold is at the
same time the $n$-fold strongly-cyclic covering of the lens space
\hbox{$L(pn-p-n,p-1)$,} branched over the $(1,1)$-knot
\hbox{$K(1,p-2,(p-1)(n-2),1)$} and the $(n-1)$-fold
strongly-cyclic covering of the lens space
\hbox{$L(p(pn-p-n),p-1)$,} branched over the $(1,1)$-knot
\hbox{$K(1,p-2,(p-1)(np-p-2),1)$,} when $p>2$. Moreover, the
Neuwirth manifold $M_n=\Sigma(n,2,1,1)=\Sigma(n-1,2,1,2)$ is at
the same time the $n$-fold strongly-cyclic covering of the lens
space \hbox{$L(n-2,1)$,} branched over the $(1,1)$-knot
\hbox{$K(1,n-2,0,1)$} and the $(n-1)$-fold strongly-cyclic
covering of the lens space \hbox{$L(2n-4,1)$,} branched over the
$(1,1)$-knot \hbox{$K(1,2n-4,0,1)$.}

\bigskip

\noindent{\bf Acknowledgements}

\noindent Work performed under the auspices of the G.N.S.A.G.A. of
I.N.d.A.M. (Italy) and the University of Bologna, funds for
selected research topics.

\vspace{15 pt} {LUIGI GRASSELLI, Department of Sciences and
Methods for Engineering, University of Modena and Reggio Emilia,
42100 Reggio Emilia, ITALY. E-mail: grasselli.luigi@unimore.it}

\vspace{15 pt} {MICHELE MULAZZANI, Department of Mathematics,
University of Bologna, I-40127 Bologna, ITALY, and C.I.R.A.M.,
Bologna, ITALY. E-mail: mulazza@dm.unibo.it}


\begin{thebibliography}{5}


\bibitem {BKM}
P. Bandieri, A. C. Kim and M. Mulazzani: {\it On the cyclic
coverings of the knot $5_2$}. Proc. Edinb. Math. Soc. {\bf 42}
(1999), 575--587.

\bibitem {BZ}
M. Boileau and H. Zieschang: {\it Heegaard genus of closed
orientable Seifert $3$-manifolds}. Invent. Math. {\bf 76} (1984),
455--468.

\bibitem {BuZ}
G. Burde and H. Zieschang: {\it Knots.} De Gruyter Studies in
Mathemathics, 5. Walter de Gruyter, 1985.

\bibitem {CJ}
A. Casson and D. Jungreis: {\it Convergence groups and Seifert
fibered $3$-manifolds}. Invent. Math. {\bf 118} (1994), 441--456.

\bibitem {CM1}
A. Cattabriga and M. Mulazzani: {\it Strongly-cyclic branched
coverings of $(1,1)$-knots and cyclic presentation of groups}.
Math. Proc. Camb. Philos. Soc. {\bf 135} (2003), 137--146.

\bibitem {CM2}
A. Cattabriga and M. Mulazzani: {\it $(1,1)$-knots via the mapping
class group of the twice punctured torus}. Adv. Geom. {\bf 4}
(2004), 263--277.

\bibitem {CM3}
A. Cattabriga and M. Mulazzani: {\it All strongly-cyclic branched
coverings of $(1,1)$-knots are Dunwoody manifolds}. J. London
Math. Soc. {\bf 70} (2004), 512--528.

\bibitem {CM4}
A. Cattabriga and M. Mulazzani: {\it Representations of
$(1,1)$-knots}. to appear in Fundamenta Mathematicae, Proceedings
of the Second International Conference ``Knots in Poland 2003'',
Vol. II (2005). arXiv:math.GT/0501234

\bibitem {Ca}
A. Cavicchioli: {\it Neuwirth manifolds and colourings of graphs}.
Aequationes Math. {\bf 44} (1992), 168--187.

\bibitem {CHK}
A. Cavicchioli, F. Hegenbarth and A. C. Kim: {\it A geometric
study of Sieradsky groups}. Algebra Colloq. {\bf 5} (1998),
203--217.

\bibitem {CHK2}
A. Cavicchioli, F. Hegenbarth and A. C. Kim: {\it On cyclic
branched coverings of torus knots}. J. Geom. {\bf 64} (1999),
55--66.

\bibitem {CHR}
A. Cavicchioli, F. Hegenbarth and D. Repovs: {\it On manifold
spines and cyclic presentations of groups}. In: Knot theory.
Proceedings of the mini-semester, Warsaw, Poland, July 13--August
17, 1995. Warszawa: Polish Academy of Sciences, Institute of
Mathematics, Banach Cent. Publ. {\bf 42} (1998), 49--56.

\bibitem {CRS}
A. Cavicchioli, D. Repovs and F. Spaggiari: {\it Topological
properties of cyclically presented groups}. J. Knot Theory
Ramifications {\bf 12} (2003), 243--268.

\bibitem {Du}
M. J. Dunwoody: {\it Cyclic presentations and 3-manifolds}. In:
Proc. Inter. Conf., Groups-Korea '94, Walter de Gruyter,
Berlin-New York (1995), 47--55.

\bibitem {GM}
L. Grasselli and M. Mulazzani: {\it Genus one 1-bridge knots and
Dunwoody manifolds}. Forum Math.  {\bf 13} (2001), 379--397.

\bibitem {GP}
L. Grasselli and S. Piccarreta: {\it Crystallizations of
generalized Neuwirth manifolds}. Forum Math. {\bf 9} (1997),
669--685.

\bibitem {HKM2}
H. Helling, A. C. Kim and J. L. Mennicke: {\it A geometric study
of Fibonacci groups}. J. Lie Theory {\bf 8} (1998), 1-23.

\bibitem {Jo}
D. L. Johnson: {\it Topics in the theory of group presentations}.
London Math. Soc. Lect. Note Ser., vol. 42, Cambridge Univ. Press,
Cambridge, U.K., 1980.

\bibitem {Ka}
A. Kawauchi: {\it A survey of knot theory.} Birkhauser Verlag,
Basel, 1996.

\bibitem {Ki}
A. C. Kim: {\it On the Fibonacci group and related topics}.
Contemp. Math. {\bf 184} (1995), 231--235.

\bibitem {KKV1}
A. C. Kim, Y. Kim and A. Vesnin: {\it On a class of cyclically
presented groups}. In: Proc. Inter. Conf., Groups-Korea '98,
Walter de Gruyter, Berlin-New York (2000), 211--220.

\bibitem {KKV2}
G. Kim, Y. Kim and A. Vesnin: {\it The knot $5_2$ and cyclically
presented groups}. J. Korean Math. Soc. {\bf 35} (1998), 961--980.

\bibitem {MR}
C. Maclachlan and A. W. Reid: {\it Generalised Fibonacci
manifolds}. Transform. Groups {\bf 2} (1997), 165--182.

\bibitem {Mu}
M. Mulazzani: {\it Cyclic presentations of groups and cyclic
branched coverings of $(1,1)$-knots}. Bull. Korean Math. Soc. {\bf
40} (2003), 101--108.

\bibitem {Ne}
L. Neuwirth: {\it An algorithm for the construction of 3-manifolds
from 2-complexes}. Proc. Camb. Philos. Soc. {\bf 64} (1968),
603--613.

\bibitem {Or}
P. Orlik: {\it Seifert manifolds}. Lecture Notes in Mathematics,
Vol. 291. Springer-Verlag, Berlin-New York, 1972.

\bibitem {RSV}
B. Ruini, F. Spaggiari and A. Vesnin: {\it On spines of Seifert
fibered manifolds}. Aequationes Math. {\bf 65} (2003), 40--60.

\bibitem{ST}
H. Seifert and W. Threlfall: {\it Seifert and Threlfall: a
textbook of topology}. Pure and Applied Mathematics, 89. Academic
Press, New York-London, 1980.

\bibitem{Si} J. Singer: {\it Three-dimensional
manifolds and their Heegaard diagrams}. Trans. Amer. Math. Soc.
{\bf 35} (1933), 88-111.

\bibitem{Sp} F. Spaggiari: {\it The combinatorics of some Tetrahedron manifolds}.
Discrete Math., to appear.

\bibitem {SV}
A. Szczepa\'nski and A. Vesnin: {\it Generalized Neuwirth groups
and Seifert fibered manifolds}. Algebra Colloq. {\bf 7} (2000),
295--303.

\bibitem {Ta}
M. Takahashi: {\it Two knots with the same $2$-fold branched
covering space.} Yokohama Math. J. {\bf 25} (1977), 91-99.

\bibitem {VK}
A. Vesnin and A. C. Kim: {\it The fractional Fibonacci groups and
manifolds}. Sib. Math. J. {\bf 39} (1998), 655--664.



\end{thebibliography}
\end{document}